\documentclass[10pt,reqno]{amsart}

\usepackage{amsmath,amsfonts,amssymb,pstricks,booktabs, enumerate,hyperref,multicol,graphicx,subfig,float,caption,subfig}

\numberwithin{equation}{section}

\begin{document}
\title[]{Symmetry analysis and soliton solutions of (2+1)- dimensional Zoomeron Equation }
\author{}
\maketitle
\begin{center}
Vishakha Jadaun,   Sachin Kumar, Yogeeta Garg
\end{center}

\begin{abstract}
 Traveling wave solutions of $(2+1)$-dimensional Zoomeron equation(ZE) are developed in terms of exponential functions involving free parameters. It is shown that the novel Lie group of transformations method is a competent and prominent tool in solving nonlinear partial differential equations(PDEs) in mathematical physics. The similarity transformation method(STM) is applied first on $(2+1)$-dimensional ZE to find the infinitesimal generators. Discussing the different cases on these infinitesimal generators, STM reduce $(2+1)$-dimensional ZE into $(1+1)$-dimensional PDEs, later it reduces these PDEs into various ordinary differential equations(ODEs) and help to find exact solutions of $(2+1)$-dimensional ZE.

\end{abstract}
 \textbf{keywords:} Zoomeron equation, Lie point transformation method, Exact solutions. 

\section{Introduction}
 Non-linear partial differential equations(PDEs) exhibit a rich variety of non-linear phenomena, arise in many physical fields like the stratified shear flow in ocean and atmosphere,condense matter physics, fluid mechanics, etc. Therefore, seeking exact solutions of  non-linear PDEs specifically the non-linear evolution equations (NLEEs) plays an important role to look into the internal mechanism of convoluted physical phenomena. Most of the physical phenomena such as, fluid mechanics, quantum mechanics, electricity, plasma physics, chemical kinematics, propagation of shallow water waves, and optical fibers are modelled by non-linear evolution equation and the appearance of solitary wave solutions in nature is somewhat frequent. The non-linear processes are one of the major challenges and not easy to control because the non-linear characteristic of the system abruptly changes due to some small changes in valid parameters including time. The phenomena of solitons was first discovered by naval engineer John Scott Russell\cite{JRS} in 1834 as a particle like wave with constant speed and shape. Later, in 1976, Calogero and Degasperis investigated the generalisation of well known KdV and Schr\"{o}dinger (1D) equation to inherit solitons and find two type of solitons Boomeron and Trappon(Zoomeron) that move with variable speed and observed a relation between their speed and polarisation effect.
 In the literature, many significant methods have been proposed for obtaining exact solutions of non-linear partial differential equations (PDEs) such as the Exp-function method, the Jacobi-elliptic method, the Lie B\"acklund transformations, the rational sine cosine method, the lie group of transformation method, the Hirota's method, Hirota bilinear forms, the tanh-sech method and so on \cite{FA, FA2,GJ1,GJ7}. The main aim of this paper is to apply the similarity transformation method to obtain soliton solutions of the (2+1)-dimensional Zoomeron(or calapso) equation (ZE) with dependent variable $u$ and independent variables $(x,y,t)$. 
 \begin{align}{\label{eq1}}
\bigg(\frac{ u(x,y,t)_{xy}}{u(x,y,t)}\bigg)_{tt} - \bigg(\frac{ u(x,y,t)_{xy}}{u(x,y,t)}\bigg)_{xx} + 2 (u^2(x,y,t))_{xt}=0
 \end{align}
 In many studies,  the forms of the arbitrary parameters which appear in the underlying model are assumed. However, the Lie symmetry approach through the method of group classification has proven to be a powerful tool in specifying the forms of these parameters naturally. The Lie group method is powerful technique to construct the exact solution of non-linear PDEs.Furthermore,based on the Lie group method, many types of exact solutions of PDEs can be considered, such as the traveling wave solutions, similarity solutions ,soliton solutions, fundamental solutions, and so on.

\section{Lie Symmetry for the Zoomeron equation}
 In this section, we perform the application of symmetry group for (2+1)-dimensional Zoomeron equation \ref{eq1} that can be written as in following form also.
 \begin{equation*}
\begin{split}
\Delta_\emph{S}(\textit{x},u^n)&=4 u^5 u_{xt}+4 u^4 u_x u_t + u^3 (u_{xytt}-u_{xxxy})+u^2(2 u_x u_{xxy}+u_{xx} u_{xy}-2 u_t u_{xyt}-u_{tt} u_{xy})\\&+ 2 u u_t^2 u_{xy}-2 u u_x^2 u_{xy}=0 , \,\,\,\,\,\,\ \emph{S}=1,.....,\emph{l}
\end{split}
\end{equation*}
Here, $\textit{x}=(x,y,t)$ and $u^n$ express all the derivatives of $u$ of all orders from $0$ to $n$. we consider One parameter Lie group of infinitesimal transformation acting on dependent variable $u$ and independent variables $(x,y,t)$ of (\ref{eq1}).
\begin{equation}
\left.
\begin{split}
x^* &=& x+\varepsilon \emph{X}(x,y,t,u)+O(\varepsilon^2)\\
y^ *&=& y+\varepsilon \emph{Y}(x,y,t,u)+O(\varepsilon^2)\\
t^* &=& t+\varepsilon \emph{T}(x,y,t,u)+O(\varepsilon^2)\\
u^* &=& u+\varepsilon \emph{U}(x,y,t,u)+O(\varepsilon^2)
\end{split}
\right\}
\end{equation}
  where $X, Y, T $ and $U$ are infinitesimals of transformation for the independent and dependent variables, respectevely and  $(\varepsilon)$ is  group  parameter which  is  admitted by  the  system  (\ref{eq1}). In other  words, Lie  group  of transformations are such  that  if $ u $ is  a  solution of equation (\ref{eq1}) then $ u^*$ is  also a solution. The method for finding group symmetry is by finding corresponding infinitesimal generator of Lie group of  transformations \cite{LVO8}. This yields to overdetermined, linear system of equation for infinitesimals $\emph{X}(x,y,t,u)$ , $\emph{Y}(x,y,t,u)$, $\emph{T}(x,y,t,u)$ and $ \emph{U}(x,y,t,u)$. The invariance of eq.\ref{eq1} under the infinitesimal transformation leads to the invariance conditions 
 \begin{equation}
 Pr^{(n)}\Delta_\emph{S}(\textit{x},u^n) =0,\emph{S}=1,.....,\emph{l} \,\,\,\,\,\ \text{whenever} \Delta_\emph{S}(\textit{x},u^n)=0
 \end{equation}
 where $Pr^{(n)}$ is called $n^{th}$-order prolongation vector field is given by
 \begin{equation}
 Pr^{(n)}=\emph{S}+\sum_{\alpha=1}^q \sum \phi_\alpha^\textit{j}(\textit{x},u^n) \partial u_\textit{j}^\alpha
\end{equation}
where $\textit{j}=(\textit{j}_1, ....., \textit{j}_k),  1\leq\textit{j}_k \leq q,  1\leq k \leq n$  , and sum is all over the orders of \textit{j}. If \textit{j}=k, the coefficients $\phi_\alpha^\textit{j}$ of $ \partial u_\textit{j}^\alpha$ will only depend on $k^{th}$ and lower order derivatives of $u$;
\begin{equation}
\phi_\alpha^\textit{j}(\textit{x},u^n)= D_\textit{j}\bigg( \phi_\alpha - \sum_{\textit{i}=1}^q \emph{X}^\textit{i} u_\textit{i}^\alpha\bigg)+\sum_{\textit{i}=1}^q  \emph{X}^\textit{i}u_{\textit{j},\textit{i}}^\alpha
\end{equation}
where $u_\textit{i}^\alpha= \frac{\partial u^\alpha}{\partial x^\textit{i}}$,  $u_{\textit{j},\textit{i}}^\alpha= \frac{\partial u_\textit{i}^\alpha}{\partial x^\textit{i}}$ and $\emph{X}^\textit{i}=(\emph{X}, \emph{Y}, \emph{T}, \emph{U})$ for eq.\ref{eq1}.

For Zoomeron equation, we consider $n=4$, thus the Infinitesimal criteria for the invariance  of Eq. \ref{eq1} would be:
\begin{equation}
 Pr^{(4)}\Delta_\emph{S}(\textit{x},u^4) =0 \,\,\, \text{where} \,\,\,  \emph{S} =\emph{X}\frac{\partial}{\partial x} +\emph{Y}\frac{\partial}{\partial y} +\emph{T}\frac{\partial}{\partial t} +\emph{U} \frac{\partial}{\partial u} 
\end{equation}
Further we assume that Eq. (1.1) is invariant. Thus surface invariant condition on the Eq. (1.1) provides
 \begin{equation}{\label{a}}
\begin{array}{lcl}
16u^{3} u_{xt} \emph{U} + 4u^{4} \emph{U}_{xt} + 12u^{2} \emph{U} u_{x} u_{t} + 4u^{3} \emph{U}_{x} \eta_{t}  +4u^{3} u_{x} \emph{U}_{t} +2 u \emph{U}( u_{xytt}- u_{xxxy}) + u^{2} (\emph{U}_{xytt}-\emph{U}_{xxxy}) + \\ \emph{U}(2 u_{x} u_{xxy} +u_{xx} u_{xy} -2 u_{t} u_{xyt}- u_{tt} u_{xy}) + u( 2 \emph{U}_{x} u_{xxy} + 2 u \emph{U}_{xxy} + \emph{U}_{xx} u_{xy} +u_{xx} \emph{U}_{xy} - \\2 \emph{U}_{t} u_{xyt} - 2u_{t} \emph{U}_{xyt} - \emph{U}_{tt} u_{xy} -u_{tt} \emph{U}_{xy}) +2 (2 u_{t} \emph{U}_{t} -2 u_{x} \emph{U}_{x})u_{xy}
 + 2(u_{t}^{2} - u_{x}^{2}) \emph{U}_{xy} =0
 \end{array}
 \end{equation}
 Substitution of extensions from \cite{PO9,GJ1} in \ref{a} provides an equation in terms of the partial derivatives of $\emph{X}$, $\emph{Y}$,$\emph{T}$ and $\emph{U}$. So one can obtain a large overdetermined system of coupled PDE’s and then we equate to zero the coefficients of various monomials. The following infinitesimals can be found after performing rather tedious calculations which can be done by hand or using any symbolic program such as Mathematica or Maple. Condition on infinitesimals can be found by solving {``determining equation"} which yields the following:
\begin{eqnarray}{\label{eq2}}
\begin{split}
\emph{U} =& -(c_{1}+c_{3}) u \\
\emph{X} =&\,\,\,\,\,2 c_{3}x+c_{5} \\
\emph{Y} =&\,\,\,\,\, 2c_{1} y+4 c_{2}\\
\emph{T} =&\,\,\,\,\, 2c_{3}t +c_{4}t\\
\end{split}
\end{eqnarray}
where $c_{1}$, $c_{2}$, $c_{3}$, $c_{4}$ and $c_{5}$ are arbitrary constants. % $\gamma(t)$ and $\lambda(t)$ are arbitrary functions of t. The prime $(\prime)$ denotes the differentiation  with respect to its indicated variable throughout the paper.\\

 The infinitesimal generators of the corresponding Lie algebra are given by
\begin{eqnarray}{\label{eq3}}
\begin{split}
V_{1} =&\,\,\, 2 y \frac{\partial}{\partial y}-u\frac{\partial}{\partial u}\\
V_{2} =&\,\,\,  \frac{\partial}{\partial y}\\
V_{3} =&\,\,\, 2 x  \frac{\partial}{\partial x}+ 2 t \frac{\partial}{\partial t}-u  \frac{\partial}{\partial u}\\
V_{4} =&\,\,\, \frac{\partial}{\partial t}\\
V_{5} =&\,\,\, \frac{\partial}{\partial x}
\end{split}
\end{eqnarray}

It is convenient to display the commutators of a Lie algebra through its commutator table whose $(i,j)^{th} $ entry is $[V_{i}, V_{j}]$.The commutator table is antisymmetric with its diagonal elements all zero as we have $[V_{\alpha}, V_{\beta}]= -[V_{\beta}, V_{\alpha}]$ (more details see \cite{GJ7, GS6}). The structure constants are easily read off from the commutator table.

For the infinitesimal generators (\ref{eq3}) we have the following commutator table: 
\begin{table}[H]
\centering
\begin{tabular}{c|ccccc}
$\cdot$ & $V_1$ & $V_2$ & $V_3$ & $V_4$ & $V_5$\\
\hline
$V_1$ & 0 & $-2V_2$ & 0 & 0 & 0\\
$V_2$ & $2V_2$ & 0 & 0 & 0 & 0\\
$V_3$ &0 & 0 & 0 & $-2V_4$ & $-2V_5$\\
$V_4$ &0&0&$2V_4$ &0&0\\
$V_5$ &0&0&0&$2V_5$ &0\\
\hline
\end{tabular}
\end{table}

Here, it is clear that the Zoomeron equation contain infinite continuous group of transformations which is generated by the infinite-dimensional Lie algebra spanned by vector fields (\ref{eq3}). These generators are linearly independent. In general, there are an infinite number of subalgebras for this Lie algebra formed from linear combinations of generators $V_{i}, i= 1,2,3,4,5$. Thus, to get the similarity solution for Eq.(\ref{eq1}), the corresponding characteristic equations are:
\begin{equation}{\label{eq4}}
\frac{dx}{\emph{X}(x,y,t,u)} = \frac{dy}{\emph{Y}(x,y,t,u)} =\frac{dt}{\emph{T}(x,y,t,u)} = \frac{dv}{\emph{U}(x,y,t,u)} 
\end{equation}
If two algebras are similar, i.e. connected to each other by a transformation from the symmetry group, then their corresponding invariant solutions are connectedd to each other by the same transformation. Therfore, it is sufficient to put all similar subalgebras into one class and select a representative from each class. The set of all these representatives is called an optimal system (For details see \cite{LVO8, PO9}). The different forms of the solution of equation $(\ref{eq1})$ are obtained from the optimal system for (\ref{eq3}) consists of following vector fields:
\begin{eqnarray}
(\mbox{i}) V_{1},\,\,\,\, (\mbox{ii}) V_{3},\,\,\,\, (\mbox{iii}) V_{4}, \,\,\,\, (\mbox{iv}) V_{2}+a V_{4}, \,\,\,\, (\mbox{v}) \alpha V_{2}-\beta V_{5}, \,\,\,\, (\mbox{vi}) V_{4}+V_{5}
\end{eqnarray}
 where $ \beta,\gamma$ and $a$ are arbitrary constants.
 \begin{table}[H]
 \label{table1}
 \caption{Similarity variables and similarity forms of Zoomeron equation.}
 \centering
 \begin{tabular}{p{2cm} p{2cm} p{2cm}}
 \hline
 Essential fields &  Similarity variables $(\mu,\delta)$ & Similarity solution $(u)$\\
 \hline
 $V_{1}$ & $x, t$ & $\frac{F(\mu, \delta)}{\sqrt{y}}$\\
 $\gamma V_{2}- \beta V_{5}$ & $\gamma x+ \beta y, y$ & $F(\mu, \delta)$\\
 $V_{3}$ & $\frac{x}{t}, y$ & $\frac{F(\mu, \delta)}{\sqrt{t}}$\\
 $V_{4}$ & $x, y$ & $F(\mu, \delta)$\\
 $V_{2}+a V_{4}$ & $x, y-\frac{t}{a}$ & $F(\mu, \delta)$\\
 $V_{4}+V_{5}$ & $x-t, y$ & $F(\mu,\delta)$\\
 [2ex]
 \hline
 \end{tabular}
 \end{table}
 The invariants and the similarity solutions of Zoomeron Equation can be find by solving Lagrange's system \ref{eq4}. The general solution of these equations involves three constants; two  are independent variables $\mu, \delta$ and the other plays the role of new independent variable$F(\mu, \delta)$.
 
\section{Reduced PDEs for different vector fields and exact solutions:}
\subsection{Vector field $V_{1}$:}
The associated Lagrange system is found by comprising (\ref{eq2}) and (\ref{eq4})
\begin{equation}
\frac{dx}{0} = \frac{dy}{2c_{1} y} =\frac{dt}{0} = \frac{du}{-c_{1}u} 
\end{equation}
The similarity reduction of equation $(\ref{eq1})$ in similarity form is
\begin{eqnarray}
u(x,y,t)=\frac{F(\mu,\delta)}{\sqrt{y}} \,\,\,
\mbox{where}\,\,\,\,\,
\mu= x \,\,\, \mbox{and} \,\,\, \delta = t 
\end{eqnarray} 
are the two invariants that we obtained.\\
From Eqs. $(2.11)$ and $(\ref{eq1})$, we get the following equation
\begin{equation}{\label{0}}
4 F^5 F_{\mu\delta} +4 F^4 F_\mu F_{\delta}+ F^3 F_{\mu\mu\mu}-\frac{F^3 F{\mu\delta\delta}}{2}- \frac{3 F^2 F_{\mu} F_{\mu\mu}}{2}+ \frac{F^2 F_{\mu} F_{\delta\delta}}{2} + F^2 F_{\delta} F_{\mu\delta}- F F_{\mu}(F_{\delta})^2+ F (F_{\mu})^3=0
\end{equation}
and the new set of infinitesimals for Eq. $(2.12)$ by applying symilarity transformation method (STM) is
\begin{eqnarray}
\xi= 2a_{1} \mu+a_{2}, \,\,\, \psi= 2 a_{1} \delta+ a_{3} , \,\,\, \eta= -a_{1}F
\end{eqnarray}
where $a_{1}$, $a_{2}$ and $a_{3}$ are arbitrary constants. This follows the characteristic equation for $(\ref{0})$ is given by
\begin{eqnarray*}
\frac{d\mu}{ 2a_{1}\mu+a_{2}}=\frac{d\delta}{2 a_{1}\delta + a_{3}}= \frac{dF}{-a_{1}F}
\end{eqnarray*}
Further F can be written as
\begin{align}
F(\mu,\delta)= \frac{H(\zeta)}{\sqrt (2 a_{1}\delta+a_{3})} \,\,\, 
 \mbox{where} \,\,\,\zeta = \frac{2a_{1}\mu+a_{2}}{2 a_{1}\delta+a_{3}} \,\,\,\, \mbox{is a similarity variable.} 
\end{align}

Then using Eq. $(\ref{0})$ and $(2.14)$, Eq.$(\ref{eq1})$ can be reduced to a third order non-linear ODE for $H(\zeta)$ as 
\begin{eqnarray}
\begin{split}{\label{*}}
-&2\zeta H^5(\zeta)  H^{\prime\prime}(\zeta) - 4H^5(\zeta)  H^{\prime}(\zeta) -2\zeta H^4(\zeta) (H^{\prime^2})(\zeta)+a_{1}(1-\frac{\zeta^2}{2}) H^3(\zeta) H^{\prime\prime\prime}(\zeta)+a_{1} H^3(\zeta) H^{\prime}(\zeta) \\ -&2 a_{1} \zeta H^3(\zeta) H^{\prime\prime}(\zeta)+2 a_{1} \zeta H^2(\zeta) H^{\prime^2}(\zeta)\frac{3}{2} a_{1}(\zeta^2-1)H^2(\zeta) H^{\prime}(\zeta)H^{\prime\prime}(\zeta)+(\zeta^2-1)H(\zeta) H^{\prime^3}(\zeta) =0
\end{split}
\end{eqnarray}
$H(\zeta)=C_{0}$ is a solution of Eq. $(\ref{*})$ and $F(\mu,\delta_{2})=\alpha(\mu)+C_{0}, \,\,\,\, F(\mu,\delta_{2})=A_{0}\mu+B_{0}+\beta(\delta_{2})$ are the solutions of Eq. $(\ref{0})$.\\ Hence, comprising Eqs. $(2.7)$, $(2.8)$ and $(2.11)$, solutions of CD equation $(\ref{eq1})$ are given by
 \begin{eqnarray}
u(x,y,t)=\frac{C_{0}}{\sqrt{y(2a_{1} t+a_{3})}}\\
\end{eqnarray}
where $C_{0}$ is an arbitrary constant.
 \subsection{Vector fields $V_{2}$ and $V_{4}$:}
 The invariants of $V_{2}$ and $V_{4}$ are respectively given by 
 \begin{align}{\label{y}}
 \mu =& x, \quad \delta = t, \quad u(x,y,t)=  F(\mu ,\delta )\quad \mbox{and} \\ 
 \mu =& x, \quad \delta = y, \quad u(x,y,t)=  F(\mu ,\delta )
 \end{align}
 By using transformation of the form $w = x+t $, \,\,\, $F(w) = u(x,y,t)$ and $w = x+y $, \,\,\, $F(w) = u(x,y,t)$ for the invariants of vector fields  $V_{2}$ and $V_{4}$ respectively, one can reduce equation (\ref{eq1}) into non-linear ODEs and solutions of (\ref{eq1}) in the form of error functions and airy funtions are as follows:
\begin{align*}
u(x,y,t)&= y^{-\frac{1}{2}} \bigg\lbrace \sqrt{\frac{2 \pi}{A}} \exp \left[\frac{-2 (m_1(x+t)+m_2)^2}{A} \right]^2\times erfi \left[\frac{1}{A} (m_1 (x+t)+ m_2)\right]\\&+ m_1 \exp \left[\frac{-2 m_1 (x+t) (m_1(x+t)+ 2 m_2)^2}{A} \right] \bigg\rbrace^{-\frac{1}{2}}\\
u(x,y,t)&= m_3 A_i \bigg\lbrace m_1^{-\frac{2}{3}}[ m_2+m_1(x+y)]\bigg\rbrace+ m_4 B_i \left[ m_2+m_1(x+y)\right]
\end{align*} 
where $m_1, m_2$ and $m_3$ are constants. $A_i$ and $B_i$ are airy functions defined by
\begin{align*}
A_i(z)&= \frac{3^{-\frac{2}{3}}}{\Gamma(\frac{2}{3})} 0F1 \left(\frac{2}{3};\frac{1}{9};z^3\right)-\frac{3^{-\frac{1}{3}}}{\Gamma(\frac{1}{3})} 0F1 \left(\frac{4}{3};\frac{1}{9};z^3\right) \quad \mbox{and}\\
B_i(z)&= \frac{3^{-\frac{1}{6}}}{\Gamma(\frac{2}{3})} 0F1 \left(\frac{2}{3};\frac{1}{9};z^3\right)-\frac{3^{-\frac{1}{6}}}{\Gamma(\frac{1}{3})} 0F1 \left(\frac{4}{3};\frac{1}{9};z^3\right)
\end{align*}
where $0F_1(a; b; z)$ is the confluent hypergeometric limit function.

 Invariants of vector field $V_{4}$ which are given in (\ref{y}), using these invariants, we can reduce (\ref{eq1}) into following PDE
 \begin{equation}{\label{eq9}}
-F^3 F_{\mu\mu\mu\delta} +2 F^2 F_{\mu}F_{\mu\mu\delta} +F^2 F_{\mu\mu}F_{\mu\delta}-2 F F_{\mu}^2 F_{\mu\delta}=0
\end{equation}
The given set of infinitesimals can be evaluated after imposing STM on (\ref{eq9})
\begin{eqnarray}{\label{new3}}
\xi=  b_{1}\mu+b_{2} , \,\,\,\, \psi= f(\delta), \,\,\,\, \eta= b_{3} F 
\end{eqnarray}
where $b_{1}$,$b_{2}$ and $b_{3}$ are arbirary constants an $f$ is the function of $\delta$
One can get the following transformation with similarity variable $\zeta$ :
\begin{eqnarray}{\label{New4}}
\begin{split}
F(\mu,\delta) =& \,\, e^\frac{-b_{3} \delta}{f(\delta)}H(\zeta)\,\,\,\,\,
\mbox{with invariant} \,\,\,\, \zeta =& ( b_{1}\mu+b_{2}) e^\frac{-b_{1} \delta}{f(\delta)}
\end{split}
\end{eqnarray}
 Eq.(\ref{eq9}) can be reduce to the following fourth order non-linear ODE
\begin{equation}{\label{10}}
\begin{split}
b_{1}\zeta H^2(\zeta)  H^{\prime\prime\prime\prime}(\zeta) -(b_{3}-3b_{1})H^2(\zeta)  H^{\prime\prime\prime}(\zeta)-2 (b_{3}-b_{1}) H^{\prime^3}(\zeta)-b_{1} \zeta H(\zeta) H^{\prime\prime^2}(\zeta) \\+3 (b_{3}-b_{1}) H(\zeta) H^{\prime}(\zeta)H^{\prime\prime}(\zeta)-2 b_{1}\zeta H(\zeta)H^{\prime}(\zeta) H^{\prime\prime\prime}(\zeta)+2 b_{1} \zeta H^{\prime^2}(\zeta)H^{\prime\prime}(\zeta)=0
\end{split}
\end{equation}
with the help of (\ref{10}), we find one more solution of Eq.(\ref{eq1})with an arbitrary constant $a$
\begin{eqnarray}
u(x,y,t)= a e^\frac{-b_{3}y}{f(y)}
\end{eqnarray} 
 
\subsection{Vector field $V_{3}$:}

 Symmetry reduction of Eq.$(\ref{eq1})$ under one point symmetry group is given by using standard method.

Firstly, solving the corresponding characteristic equations by comprising Equations (\ref{eq2}) and (\ref{eq4})
\begin{equation}
\frac{dx}{2c_{3} x} = \frac{dy}{0} =\frac{dt}{2 c_{3}t} = \frac{dv}{-c_{3} u} 
\end{equation}
we obtain two invariants with given similarity solution:
\begin{eqnarray}
\begin{split}{\label{eq5}}
 \mu=& \frac{x}{t} \,\,\,\,\, \mbox{and} \,\,\,\,\, \delta = y\\
 u(x,y,t)=& \frac{ F(\mu,\delta)}{\sqrt{t}} 
 \end{split}
 \end{eqnarray}

 Eq.$(\ref{eq1})$ can be reduce to the following PDE, using Eq.(\ref{eq5})
\begin{equation}{\label{eq6}}
\begin{split}
-&4 \mu F_{\mu\mu} F^5 -8 F_{\mu} F^5-4\mu F_{\mu}^2 F^4+2 F^3 F_{\mu\delta}+ (\mu^2-1)F^3 F_{\mu\mu\mu\delta}+ 4 \mu F^3 F{\mu\mu\delta}\\+& 2 (1-\mu^2) F^2 F_{\mu} F_{\mu\mu\delta}+(1-\mu^2) F^2 F_{\mu\mu} F_{\mu\delta}-4 \mu F^2 F_{\mu} F_{\mu\delta} - 2 (1-\mu^2) F F_{\mu}^2 F_{\mu\delta}=0
\end{split}
\end{equation}
Further, applying STM, one can get a new set of infinitesimals
\begin{eqnarray}{\label{eq7}}
\xi=  0 , \,\,\,\, \psi=\overset{*}{a_{1}} \delta_{2}+ \overset{*}{a_{2}}, \,\,\,\, \eta= -\overset{*}{a_{1}} F 
\end{eqnarray}
where $\overset{*}{a_{1}}$ and  $\overset{*}{a_{2}}$, a new set of parameters, obtained by further applying STM.

Solve the following Lagrange's system corresponding to infinitesimals (\ref{eq7})
\begin{align*}
\frac{d\mu}{0} = \frac{\delta}{\overset{*}{a_{1}} \delta+ \overset{*}{a_{2}}} = \frac{dF}{-\overset{*}{a_{1}} F}
\end{align*}
It gives the transformation of the form
\begin{eqnarray}{\label{New1}}
\begin{split}
F(\mu,\delta) =& \frac{H(\zeta)}{\sqrt{(\overset{*}{a_{1}} \delta+ \overset{*}{a_{2}})}}\,\,\,\,\,\,\,
\mbox{with invariant} \,\,\,\, \zeta =& \mu
\end{split}
\end{eqnarray}
which reduces Eq.(\ref{eq6}) to the following third order non-linear ODE
\begin{equation}{\label{New2}}
\begin{matrix}
\overset{*}{a_{1}}(\zeta^2-1) H^2(\zeta)  H^{\prime\prime\prime}(\zeta) - 4 \overset{*}{a_{1}}\zeta H^2(\zeta) H^{\prime\prime}(\zeta) +4 \zeta H^3(\zeta) H^{\prime^2}(\zeta)&\\ +3 \overset{*}{a_{1}} (\zeta^2-1) H(\zeta) H^{\prime}(\zeta)H^{\prime\prime}(\zeta)+4 \zeta H^4(\zeta) H^{\prime\prime}(\zeta)+8 H^4(\zeta) H^{\prime}(\zeta) \\
+2 \overset{*}{a_{1}} H^2(\zeta) H^{\prime}(\zeta)-4 \overset{*}{a_{1}}\zeta H(\zeta) H^{\prime^2}(\zeta)+2 \overset{*}{a_{1}}(\zeta^2-1)H^{\prime^3}(\zeta)  
\end{matrix}=0
\end{equation}
 solutions of (\ref{New2}) with constant $C_{3}$ are following 
\begin{align*}
 H(\zeta)=& \,\, C_{3}\\
 H(\zeta)=&\,\, \sqrt{\frac{a_{1}}{6\zeta_{1}^2}}
 \end{align*}
Hence,we find the solutions of $(\ref{eq1})$ which are given by
\begin{eqnarray}
u(x,y,t)&=& \frac{C_{3}}{\sqrt{t(\overset{*}{a_{1}} y+ \overset{*}{a_{2}})}} \\
u(x,y,t) &=&  \sqrt{\frac{a_{1}t}{6 x^2(\overset{*}{a_{1}} y+ \overset{*}{a_{2}})}}
\end{eqnarray}
\subsection{Vector field $V_{2}+a V_{4}$} For subalgebra $V_{2}+a V_{4}$, invariants are: 
\begin{eqnarray}
\mu = x,\,\,\,\,\,\ \delta = y-\frac{t}{a}, \,\,\,\,\   u(x,y,t)= F(\mu, \delta)
\end{eqnarray}
Therefore, equation \ref{eq1} reduces into following PDE:
\begin{equation}{\label{n1}}
\begin{split}
-&4 a F^4 F_{\mu\delta} -4 a F^3 F_{\mu} F_{\delta} +F^2 F_{\mu \delta \delta \delta}- a^2 F^2 F_{\mu \mu \mu \delta}+2 a^2 F F_{\mu} F_{\mu\mu\delta}\\+& a^2 F F_{\mu \mu} F{\mu\delta}- 2 F F_{\delta} F_{\mu\delta\delta}- F F_{\delta \delta} F_{\mu\delta}+2 (F_{\delta}^2 - a^2 F_{\mu}^2) F_{\mu\delta}=0
\end{split}
\end{equation}
The Lie point symmetry generators of equation \ref{n1} are given as follows:
\begin{eqnarray}{\label{n2}}
\xi= a_{1} \mu +a_{2}, \,\,\, \psi= a_{1} \delta+ a_{3} , \,\,\, \eta= -a_{1} F
\end{eqnarray}
where $\xi$, $\psi$ and $\eta$ are infinitesimals corresponding to $F$, $\mu$ and $\delta$ respectively and $a_{1}$, $a_{2}$ and $a_{3}$  are arbitrary constants.
Here we utilize the Lie point symmetry generators of equation \ref{n1} found in \ref{n2} to obtain symmetry reduction with group invariant $\zeta$ 
 \begin{eqnarray}
F(\mu,\delta) = \frac{H(\zeta)}{(a_{1} \delta+ a_{3})} \,\,\,\,\,\ \mbox{with invariant,} \,\,\,\,\,
 \zeta = \frac{a_{1} \mu + a_{2}}{a_{1} \delta + a_{3}}
\end{eqnarray}
and construct exact group-invariant solutions for equation \ref{n1} by reducing equation \ref{n1} to an ODE.
\begin{equation}{\label{n3}}
\begin{split}
&a_{1}^2 (\zeta-\zeta^3)H^2(\zeta) H^{\prime\prime\prime\prime}(\zeta)+a_{1}^2(4-10 \zeta^2)H^2(\zeta) H^{\prime\prime\prime}(\zeta)-a_{1}^2(\zeta-\zeta^3) H(\zeta) H^{\prime}(\zeta) H^{\prime \prime\prime}(\zeta)\\+&4\zeta H^2(\zeta) H^{\prime\prime}(\zeta) +12 H^4(\zeta)  H^{\prime}(\zeta)+4 \zeta H^3(\zeta)-26 a_{1}^2 H^2(\zeta) H^{\prime}(\zeta)-4 a_{1}^2 (\zeta^2-1)H^{\prime^3}(\zeta)\\-&22 a_{1}^2 \zeta H^2(\zeta)H^{\prime\prime}(\zeta)+a_{1}^2 (18\zeta^2-8)H(\zeta) H^{\prime}(\zeta)H^{\prime\prime}(\zeta)+\zeta^2 H(\zeta) H^{\prime\prime^2}(\zeta)=0
\end{split}
\end{equation}

\subsubsection{Vector field $ \gamma$ $V_{2}-$ $\beta$ $V_{5}$:} In this case we get the following Lagrange system
\begin{equation}
\frac{dx}{-\beta} = \frac{dy}{\alpha} =\frac{dt}{0} = \frac{du}{0} 
\end{equation}
The invariants for this subalgebra has already been defined in Table 1. Therefore, the reduced PDE in this case is given by
\begin{equation}{\label{7}}
\begin{split}
&-\beta F^2 F_{\mu\mu\mu \mu} - F_{\mu \mu \mu \delta} F^2+2 \beta F F_{\mu} F_{\mu \mu \mu} + 2 F F_{\mu} F_{\mu\mu\delta}\\ &+ \beta F F{\mu\mu}^2+ F F_{\mu \delta} F_{\mu\mu}-2 \beta F_{\mu}^2 F_{\mu\mu}-2 F_{\mu}^2 F_{\mu\delta}=0
\end{split}
\end{equation}
Apply Lie Symmetry method to the equation \ref{7} to get new set of infinitesimals as follows:
\begin{eqnarray}
\xi= \frac{-a_{1}}{\beta} \mu+ a_{1} \delta +\beta P(\delta) + a_{2}, \,\,\, \psi= P(\delta) , \,\,\, \eta= 0
\end{eqnarray}
where $\xi$, $\psi$ and $\eta$ are infinitesimals corresponding to $F$, $\mu$ and $\delta$ respectively and $a_{1}$, $a_{2}$ are arbitrary constants and $P(\delta)$ is an arbitrary function of $\delta$.
\subsubsection{$P(\delta)= C_{0}$}
 To obtain the symmetry reductions of equation \ref{7}, we have to solve the following characteristic equation:
 \begin{align*}
\frac{d\mu}{\frac{-a_{1}}{\beta} \mu+ a_{1} \delta +\beta C_{0}} = \frac{\delta}{C_{0}} = \frac{dF}{0}
\end{align*}
We use the method of characteristics to determine the symmetry variables and the invariants:
\begin{eqnarray}{\label{0}}
F(\mu,\delta) = H(\zeta) \,\,\,\,\,\,\  \mbox{with invariant,} \,\,\,\,\,
 \zeta = e^{\frac{a_{1} \delta}{C_{0} \beta}} \bigg(\mu-\beta \delta - \frac{a_{2} \beta}{a_{1}}\bigg) 
\end{eqnarray}
Now, reduction of variable is performed to obtain ODE:

\begin{equation}
\begin{split}
&\zeta H^2(\zeta) H^{\prime\prime\prime\prime}(\zeta)- 2 \zeta   H(\zeta) H^{\prime}(\zeta) H^{\prime\prime\prime}(\zeta)-4 H(\zeta) H^{\prime}(\zeta) H^{\prime\prime}(\zeta)\\+& 3 H(\zeta) H^{\prime\prime\prime}-H(\zeta) H^{\prime\prime}(\zeta)-\zeta H^{\prime\prime^2}(\zeta)+2\zeta H^{\prime^2}(\zeta) H^{\prime\prime}(\zeta)+2 H^{\prime^3} =0
\end{split}
\end{equation}

\subsection{Vector field $V_{4}+V_{5}$:}  The reduced PDE for this subalgebra is given by:
\begin{equation}{\label{8}}
F_{\mu \mu} F_{\mu \delta}-4 F^3 F_{\mu \mu}-4 F^2 F_{\mu}^2 -F_{\mu \mu}^2=0
\end{equation}
Apply Lie Symmetry method to the equation \ref{8} as already applied. In this case we get symmetries as under:
\begin{eqnarray}
\xi= (a_{1}+a_{3}) \mu+a_{3} \delta + a_{4}, \,\,\, \psi= a_{1} \delta+ a_{2} , \,\,\, \eta= -(a_{1}+\frac{a_{3}}{2}) F
\end{eqnarray}
where $\xi$, $\psi$ and $\eta$ are infinitesimals corresponding to $F$, $\mu$ and $\delta$ respectively and $a_{1}$, $a_{2}$, $a_{3}$ and $a_{4}$ are arbitrary constants. Thus the symmetry algebra admitted by \ref{8} is
\begin{eqnarray}{\label{8.1}}
\begin{split}
\overset{*}{V_{1}} =&\,\,\, \mu \frac{\partial}{\partial \mu}+ \delta \frac{\partial}{\partial \delta}-F \frac{\partial}{\partial F}\\
\overset{*}{V_{2}} =&\,\,\,  \frac{\partial}{\partial \delta}\\
\overset{*}{V_{3}} =&\,\,\, (\mu + \delta) \frac{\partial}{\partial \mu}- \frac{F}{2} \frac{\partial}{\partial F}\\
\overset{*}{V_{4}} =&\,\,\, \frac{\partial}{\partial \mu}
\end{split}
\end{eqnarray}
The commutator table in this case will be:
\begin{table}[H]
\centering
\begin{tabular}{c|cccc}
$\cdot$ &$\overset{*}{V_{1}}$ & $\overset{*}{V_{2}}$ & $\overset{*}{V_{3}}$ & $\overset{*}{V_{4}}$ \\
\hline
$\overset{*}{V_{1}}$ & 0 & $-\overset{*}{V_{2}}$ & 0 & $-\overset{*}{V_{4}}$ \\
$\overset{*}{V_{2}}$ & $\overset{*}{V_{2}}$ & 0 & $\overset{*}{V_{4}}$ & 0 \\
$\overset{*}{V_{3}}$ &0 & $-\overset{*}{V_{4}}$ & 0 & $-\overset{*}{V_{4}}$ \\
$\overset{*}{V_{4}}$ &$\overset{*}{V_{4}}$ &0 &$\overset{*}{V_{4}}$ &0\\
\hline
\end{tabular}
\end{table}

In general, there are infinite number of subalgebras of this Lie algebra formed from any linear combination of generators $Vj$,  $j = 1, 2, 3, 4$ and to each sublagebra one can get the reduction using characteristic equations:
\begin{align*}
\frac{d\mu}{(a_{1}+a_{3}) \mu+a_{3} \delta + a_{4}} = \frac{\delta}{a_{1} \delta+ a_{2}} = \frac{dF}{-(a_{1}+\frac{a_{3}}{2}) F}
\end{align*}
We again use the method of characteristic  to determine the following invariants and the reduced ODE.
 \begin{eqnarray}{\label{n}}
\begin{split}
F(\mu,\delta) &= \frac{H(\zeta)}{(a_{1} \delta+ a_{2})^\frac{2a_{1}+a_{3}}{2 a_{1}}}\\ \mbox{with invariant,} \,\,\,\,\,
 \zeta &= \frac{a_{1} \mu}{(a_{1} \delta+ a_{2})^\frac{a_{1}+a_{3}}{ a_{1}}} + \frac{1}{(a_{1} \delta+ a_{2})^\frac{a_{3}}{ a_{1}}} +\frac{Q}{(a_{1} \delta+ a_{2})^\frac{a_{1}+a_{3}}{ a_{1}}}
\end{split}
\end{eqnarray}
where $Q= \frac{a_{1} a_{4}-a_{2} a_{3}}{a_{1}+a_{3}}$. Substituting \ref{n} into equation \ref{8}, we obtain an ODE in $\zeta$ given as:
\begin{equation}
a_{1} (4 a_{1}+3 a_{3})H^{\prime}(\zeta)H^{\prime\prime}(\zeta)+ a_{1} (a_{1}+ a_{3}) \zeta H^{\prime\prime^2}(\zeta)+4 H^3(\zeta)H^{\prime\prime}(\zeta)+4 H^2(\zeta) H^{\prime^2}(\zeta)=0
\end{equation}

\section{Conclusion} In this paper we have shown that the Zoomeron equation can be transformed by a point transformation to fourth/third order non-linear ODEs. The Lie point symmetry generators of the ZE  were obtained by using the Lie symmetry group analysis.The analytical properties of the solutions such as  travelling wave and solitons are discussed . This work is significant since the exact solutions so obtained shall be helpful  in other applied sciences as condensed matter physics, field theory, fluid dynamics, plasma physics, non-linear optics, etc. where solitons and periodic structures are involved. Our exact solutions may serve as benchmark in the accuracy testing and comparison of their numerical algorithms. The availability of computer systems like Mathematica or Maple facilitates the tedious algebraic calculations. The method which we have proposed in this article is also a standard, direct and computer-literate method, which allows us to solve complicated and tedious algebraic calculation.

Acknowledgments:  The first author also expresses her gratitude to the University Grants Commission, New Delhi, India for financial support to carry out the above work.

\begin{figure}[H] 
\centering
%\subfloat[the figure of travelling wave solution for Eq.(3.21)]{\includegraphics[width=4cm, height=4cm]{"ZME/1"}} 
%\qquad
%\subfloat[a figure of singularity solution of (3.22)]{\includegraphics[width=4cm, height=4cm]{"ZME/2"}}  
\end{figure}
\end{document}